\documentclass{jnmp03c}

\begin{document}
\setcounter{page}{54}

\renewcommand{\evenhead}{I~A~Shereshevskii}
\renewcommand{\oddhead}{Orthogonalization of Graded Sets of Vectors}

\thispagestyle{empty}

\FistPageHead{1}{\pageref{shereshevskii-firstpage}--\pageref{shereshevskii-lastpage}}{Letter}

\copyrightnote{2001}{I~A~Shereshevskii}

\Name{Orthogonalization of Graded Sets of Vectors}
\label{shereshevskii-firstpage}

\Author{I~A~SHRESHEVSKII}

\Adress{Institute for Physics of Microstructures, 
Russian Academy of Sciences\\
 46 Uljanova street, Nizhnii Novgorod, RU-603600, Russia\\
E-mail: ilya@ipm.sci-nnov.ru}

\Date{Received September 16, 2000; Accepted October 18, 2000}

\begin{abstract}
\noindent
I propose an orthogonalization procedure preserving the grading of the
initial graded set of linearly independent vectors.  In particular,
this procedure is applicable for orthonormalization of any countable
set of polynomials in several (finitely many) indeterminates.
\end{abstract}

\noindent There are two well-known procedures for orthogonalization of
any set of linearly independent vectors in a linear space.  The first
of them boils down, essentially, to calculation of the
$\left(-\frac12\right)$-th power of the Gram matrix of the initial set
(we refer to it as {\it Gram method}), and another one is the {\it
Gram--Schmidt process}.  The result of the Gram--Schmidt process
essentially depends of the {\it order} of the elements in the set and
the operations with Gram matrix seem to be impossible to perform in
infinite dimensional spaces.

There are, however, some situations when one has to orthogonalize a
finite or infinite set of linearly independent vectors in a linear
space in such a way that the result would have some properties of the
initial set.  Observe that the set considered is not necessarily a
basis.  Any set consisting of linearly independent vectors will do. 
Observe also that I deal here with algebraic problems, no analytical
problem (convergence, completeness, etc.)  arises.

We consider some examples:
\begin{itemize}
\topsep0mm \partopsep0mm \parsep0mm \itemsep0mm \item[1.]  Let
$\left\{e_{m}= e^{imx}\right\}_{m=-\infty}^{\infty}$ be a set of
linearly independent vectors in $L_{2}([0, 2\pi], \rho)$, where $\rho$
is any positive weight.  If $\rho$ is not constant, the elements of
this set are not orthogonal, but have an important property $e_{m}
=\overline{e_{-m}}$, where the bar means complex conjugation.  Is it
possible to preserve this property under orthogonalization?  

\item[2.] Let $x=(x_{1},\ldots , x_{n})$, where $x_{j}\in {\mathbb R}$, and
$m=(m_{1}, \ldots, m_{n})$, where $m_{j}\in {\mathbb Z}_{+}$.  Then
the set $\left\{x^{m}: m\in{\mathbb Z}_{+}^{n} \right\}$, where
$x^{m}=x_{1}^{m_{1}}\cdots x_{n}^{m_{n}}$, is, due to
Stone--Weierstrass theorem, a basis in $L_{2}(\Omega)$ for a ``good''
bounded domain $\Omega\subset{\mathbb R}^{n}$.  This basis has no
natural order, but has a natural grading: the degree of $x^{m}$ is
equal to $|m|=m_{1}+\cdots +m_{n}$.  How to orthogonalize this basis
and preserve the natural grading?  Is this possible?
\end{itemize}

It is very strange that a very simple and natural answer to these and
similar questions seems to be unknown.  For this reason I present here
the orthogonalization procedure, which combines some features of the
Gram method and Gram--Schmidt process (in particular, it can be
applied in infinite dimensional case) and gives a solution to the
above problems and similar ones.  Note in this connection that,
although the orthogonalization problem is more than hundred years old,
various aspects of the propblem arise from time to time in connection
with very interesting practical problems, see, e.g., \cite{Sriv}.

Let $\left\{e^{k}_{\alpha}: \alpha\in I_{k}, \; k\in {\mathbb
Z}_{+}\right\}$ be a set of linearly independent vectors in a Hilbert
space.  We assume that all index sets $I_{k}$ are finite.  Let us
inductively define the sets of vectors $\left\{f^{k}_{\alpha}:
\alpha\in I_{k}\right\}$ by the formula
\begin{equation}
    f^{k}_{\alpha}=\sum_{\beta\in I_{k}}Q^{k}_{\beta\alpha }e^{k}_{\beta}+
                             \sum^{k-1}_{j=0}\sum_{\beta\in
                             I_{j}}P^{kj}_{\beta\alpha }f^{j}_{\beta},
    \label{base}
\end {equation}
where the unknown matrices $Q^{k}$ and $P^{kj}$ are
determined from the orthonormality conditions for the system $f$:
\begin{equation}
\ba{l}
    \left(f^{k}_{\alpha},f^{k}_{\beta}\right)=\delta _{\alpha\beta}, \qquad
    \alpha,\beta\in I_{k},
\vspace{2mm}\\
    \left(f^{k}_{\alpha},f^{j}_{\beta}\right)=0, \qquad
    \alpha\in I_{k},\qquad \beta\in I_{j},\qquad  j\neq k.
    \ea
    \label{orth}
\end{equation}

Note, first of all, that the system $f$ is linearly independent if
and only if the matrices~$Q^{k}$ are nondegenerate.  We will use
this fact later.

From the last line in (\ref{orth}) and the definition (\ref{base})
one immediately deduces that
\begin{equation}
   P^{kj}_{\beta\alpha }=
  -\sum_{\gamma\in I_{k}}D^{kj}_{ \beta\gamma}Q^{k}_{\gamma\alpha },
        \label{pdef}
\end {equation}
where we define the matrix $D^{kj}$ to be $D^{kj}_{
\beta\gamma}=\left(e^{k}_{\gamma}, f^{j}_{\beta}\right)$. Let us
substitute this expression for~$P$ in terms of  $Q$ and $D$ into
the first line in (\ref{orth}). We obtain, after
simplification, a matrix equation of the form
\begin{equation}
Q^{k\dagger}B^{k}Q^{k}=E,
\; {\rm where}\; B^{k}=\Gamma^{k}-\sum^{k-1}_{j=0}\Delta^{kj}
\label{gram} \end {equation} 
and where $\dagger$ denotes the Hermitian conjugation,
$\Delta^{kj}=D^{kj\dagger}D^{kj}$,
$\Gamma^k_{\alpha\beta}=(e^{k}_{\alpha},e^{k}_{\beta})$ is the Gram matrix
of the system $\{e^k_{\alpha}\}_{\alpha\in I_k}$ and $E$ is the unit
matrix.

Note that the matrix $B$ is the Gram matrix for the linearly
independent system of vectors in the Hilbert space
\[
    h^{k}_{\alpha}=e^{k}_{\alpha}-
               \sum^{k-1}_{j=0}\sum_{\beta\in
               I_{j}}\left(e^{k}_{\alpha},f^{j}_{\beta}\right)f^{j}_{\beta}.
\]
Hence, $B^{k}$ is positive definite.  So one can write the unique positive
definite solution of the equation (\ref{gram}) in the ``Gram'' form
\begin{equation}
      Q^{k}=\left(\Gamma^{k}-\sum^{k-1}_{j=0}\Delta^{kj}\right)^{-1/2}.
        \label{answ}
\end {equation}
This completes the orthonormalization process.

We consider now some simple examples.
\begin{itemize}
\topsep0mm
\partopsep0mm
\parsep0mm
\itemsep0mm
\item[3.] If for all $k$ the sets $I_{k}$ are one-element sets, then the
process described is exactly the Gram--Schmidt one.
\item[4.]  In the above procedure ${\mathbb Z}_{+}$ can be replaced with
any its finite subset.  If such subset is a one-element set, then our process
is exactly the Gram method.
\end{itemize}

These two examples show that the process suggested is simply a
combination of two well-known processes.
\begin{itemize}
\topsep0mm
\partopsep0mm
\parsep0mm
\itemsep0mm
\item[5.]  Let $I_{0}=\{0\}$ and $I_{k}=\{+, -\}$ for nonzero $k$'s.  Let
further $e^{k}_{\alpha}\in L_{2}([0, 2\pi], \rho)$ and $e^{0}=1$,
$e^{k}_{\pm}=e^{\pm imx}$.  (This is Example 1.)  Then, clearly,
$f^{0}=\mbox{const}\in {\mathbb R}$.  Let us show now that if
$f^{j}_{+}=\overline{f^{j}_{-}}$ for all $j<k$, then the same is true
for $k$ also.
\end{itemize}

Indeed, due to relations (\ref {orth}) and (\ref {pdef}) we obtain
\[
\ba{l}
\ds     f^{k}_{\pm} =Q^{k}_{\pm +}\left(e^{k}_{+} - \left(e^{k}_{+},f^{0}\right)f^{0}\right)
                 -\sum_{j=1}^{k-1}\left(\left(e^{k}_{+},f^{j}_{+}\right)f^{j}_{+}+
         \left(e^{k}_{+},f^{j}_{-}\right)f^{j}_{-}\right)
\vspace{3mm}\\
\ds \phantom{f^{k}_{\pm} =} +
                 Q^{k}_{\pm -}\left(e^{k}_{-} - \left(e^{k}_{-},f^{0}\right)f^{0}\right)
                 -\sum_{j=1}^{k-1}\left(\left(e^{k}_{-},f^{j}_{+}\right)f^{j}_{+}+
         \left(e^{k}_{-},f^{j}_{-}\right)f^{j}_{-}\right).
\ea
\]
Now using inductive hypothesis and the fact that $Q$ are Hermitean
matrices, i.e., $Q^{k}_{+ +}$, $Q^{k}_{--}\in {\mathbb R}$ and
$Q^{k}_{+ -}=\overline {Q^{k}_{-+}}$, it is easy to see that the
vectors $f^{k}_{\pm} $ are complex conjugate.

So we have obtained a simple affirmative answer to the question in
Example~1.

It is interesting whether or not the above construction can
be generalized to the vector systems in pseudo-euclidean spaces?  Two
problems arise in this case:
\begin{itemize}
\topsep0mm
\partopsep0mm
\parsep0mm
\itemsep0mm
\item[1)]  What shall we do with ``isotropic'' vectors?
\item[2)]  For which matrix in the right hand side instead of the identity
one, is equation (\ref{gram}) solvable?
\end{itemize}

Note in this connection that the Gram--Schmidt orthogonalization
process is not applicable in the pseudo-euclidean case because even if
the initial vector system does not contain isotropic vectors, such
vectors can appear under execution of the process and terminate it.

We consider an example: if vector $e_{1}$ is isotropic, there does not
exist any number $\alpha$ such that vector $e_{1}+\alpha e_{2}$ is
pseudo-orthogonal to $e_{1}$ (of course, unless $e_{1}$ and $e_{2}$
are initially orthogonal).

Certain properties of bases and linear independent systems in
pseudo-euclidean spaces are discussed in \cite{B,GLR}.  In
particular, Bogn\' ar proposes a ``forceful'' method for
generalization of Gram--Schmidt process to pseudo-euclidean case for
linearly ordered systems of vectors~\cite{B}.  We use something like
his method for graded systems.

In what follows we suppose that the pseudo-euclidean scalar product is
{\it ``nondegene\-rate''}, i.e., {\it for any linearly independent
finite system of vectors consisting of more than one element its Gram
matrix is nondegenerate.} In other words, this means that the
dimension of any maximal isotropic subspace (all of whose vectors are
isotropic) does not exceed~$1$.  Then any {\it finite} system of
vectors, excluding one-element ones, may be ``pseudo-orthonormalized''
in the sense that the pseudo-norm, or ``length'', of each final vector
will be equal to~$\pm 1$.

To this end, it suffices to determine the signature $(p, q)$ of the
Gram matrix $\Gamma$ of the initial system of vectors and then solve
the equation
\be\label{eq:*}
               R^{\dagger}\Gamma R=E_{p}\oplus \left(-E_{q}\right)
\ee
for unknown matrix $R$.  If we write the Hermitean matrix $\Gamma$ in
the form $\Gamma=U\Lambda U^{\dagger}$, where $U$ is a unitary matrix
and $\Lambda=\Lambda_{p}\oplus\left(- \Lambda_{q}\right) $ is a diagonal matrix
with $p$ positive and $q$ negative elements, then the solution of
(\ref{eq:*}) can be written in the form 
$R=U\left(\Lambda_{p}^{-1/2}\oplus\left(-\Lambda_{q}^{-1/2}\right)\right)$.

Therefore, it is clear that the only trouble which might appear when
dealing with graded systems of vectors in pseudo-euclidean spaces is
when after the $k$-th step the set $\left\{f^{k}_{\alpha}:\ \alpha \in
I_{k}\right\}$ consists of exactly one isotropic vector.  There does not
exists any general natural way to resolve this situation and details
depend on the concrete case.  A most simple idea is to include such
a vector into the next, $(k+1)$-th, level set of vectors and then
continue the process.  This trick slightly violates the initial
grading, but preserves the filtration.

The problems discussed stem from several sources.  One
is numerical analysis and data processing.  Here the role of
discrete Fourier transformation is well-known, but what should we
do if our data is given on a nonuniform grid?  For the answer see
Example~1.

Another problem is orthonormalization of splines which constitute a
not linearly ordered, but a graded set of functions.

D~Leites and A~Sergeev pointed out a totally different area in which
the same question arises.  These problems concern with new
polynomials in several indeterminates connected with some Lie algebras
and superalgebras and the space of these polynomilas is naturally
endowed with a nondegenerate indefinite metric, see~\cite{LS}.  For
one indeterminate D~Leites and A~Sergeev can orthogonalize their
polynomials; for several indeterminates these polynomials are not
linearly ordered and they got stuck.  By our method one can
orthogonalize the polynomials in several indeterminates proposed
in~\cite{LS}.

For the reader who wishes to compare various orthogonalization methods
I suggest very transparent and user friendly paper by Srivastava
\cite{Sriv}.

\medskip

{\bf Acknowledgements.}
I am thankful to D~Leites and Stockholm University for
hospitality, TBSS for financial support, and to T~Ya~Azizov for
helpful comments.


\label{shereshevskii-lastpage}

\end{document}